\documentclass[a4paper,12pt]{article}

\usepackage{amsfonts}

\def\ds{\displaystyle}

\def\a{\alpha}

\def\d{\delta}
\def\s{\sigma}
\def\t{\tau}
\def\tp{{\tau_p}}

\def\r{\rho}

\def\e{\varepsilon}

\def\E{{\bf E}}
\def\Ex{\E^x}

\def\F{{\mathcal F}}

\def\R{{\mathbb R}}
\def\L{{\mathbb L}}

\def\C{{\mathcal C}}
\def\M{{\mathfrak M}}

\def\I{{\mathcal I}}
\def\S{{\mathcal S}}

\def\beq{\begin{equation}}
\def\eeq{\end{equation}}

\begin{document}

\begin{center}
\Large\bf
Threshold Strategies in Optimal Stopping Problem for Diffusion Processes and Free-Boundary Problem\footnote{The results of this paper were partially presented at EPSRC Symposium Workshop {\it Optimal stopping, optimal control and finance} (Warwick, July 16--20,  2012), International Conference {\it Stochastic Optimization and Optimal
Stopping} (Moscow, September 24--28,  2012), and The Seventh Bachelier Colloquium on Mathematical finance and Stochastic calculus (Metabief, January 13--20,  2013)}
\\[0.5cm]
\large\rm
V. I. Arkin\footnote{Central Economics and Mathematics Institute, Russian Academy of Sciences (CEMI RAS). 117418, Nakhimovskii pr-t, 47, Moscow, Russia. E-mail: arkin@cemi.rssi.ru}, A. D. Slastnikov\footnote{CEMI RAS. E-mail: slast@cemi.rssi.ru}
\end{center}


\begin{abstract}

We study a problem when a solution to optimal stopping problem for one-dimensional diffusion will generate by threshold strategy. Namely, we give necessary and sufficient conditions under which an optimal stopping time can be specified as the first time when the process exceeds some level (threshold), and a continuation set is a semi-interval. We give also second-order conditions, which allow to discard such solutions to free-boundary problem that are not the solutions to optimal stopping problem.\\[4pt]

\end{abstract}


{\bf 1. Introduction.} \quad Suppose that, on some stochastic basis $(\Omega, \F,
\{\F_t, t \ge 0\},{\bf P})$, a homogeneous diffusion process $X_t,\ t\ge 0$,
with values in the interval ${\mathcal I}=]l,r[\subset{\mathbb R}^1$, where $-\infty\leq l<r\leq\infty$, is given and is described by the stochastic differential equation
\beq
dX_t=a(X_t)dt + \s(X_t)dW_t, \quad
X_0=x, \label{stochproc}
\eeq
where $a:\ \I \to \R^1$ and $\s:\ \I \to \R^1_+$ are the drift and diffusion functions and $W_t$ is the standard Wiener process.
The process $X_t$ is assumed to be regular; this means
that, starting from an arbitrary point $x\in \I$, this process
reaches any point $y\in \I$ in finite time with positive
probability.

For example, these assumptions are guaranteed, if the drift and diffusion functions satisfy the following local integrability condition
$$
\int_{x-\e}^{x+\e}\frac{1+|a(y)|}{\s^2(y)} dy <\infty \quad \mbox{for some }\e>0,
$$
at every $x\in \I$ (see, e.g., \cite{KS}).

We also assume that the boundary point $l$ of the
process values is either a natural boundary or an entry-%
not-exit. This means that it cannot be reached from an
interior point of the range of values of the process $[l,r]$
(for more details, see, e.g., \cite[Chapter 2]{BorSal}).

The process $X_t$ defined by stochastic differential equation (\ref{stochproc}) is associated with the infinitesimal
operator
\beq\label{operator}
\L f(x)= a(x)f'(x) +\frac12\s^2(x) f''(x).
\eeq

As is known, for the ODE $\L u(x)=\r u(x)$ on the
interval $\I$, there exist (unique up to constant positive multipliers) increasing and decreasing solutions $\psi(x)$ and
$\varphi(x)$, respectively (see \cite[Chapter 2]{BorSal}). Moreover, under the above assumptions on the boundary points of the interval $\I$ we have $0<\psi(x),\ \varphi(x)<\infty$ for $x\in \I$, and $\psi(l+0)=0$.
\medskip

Consider the following optimal stopping problem
for this process:
\beq \label{optstop}
\Ex g(X_\t) e^{-\r\t}\to \sup_{\t\in \M}  ,
\eeq
where $g:\ \I \to \R^1$ is the payoff function, $\r\ge 0$ is the discount rate, and $\Ex$ is mathematical expectation for the process $X_t$ which starts at the point $x$,
and the maximum is taken over the class $\M$ of stopping times $\t$ (with respect to the natural filtration
$\F_t^X=\s \{X_s,\
0\le s\le t\},\ t\ge 0$).

The function $g(x)$ is assumed to be continuous, bounded below, and $g(x_0) > 0$ for some $x_0\in \I$.\\

{\bf 2. Background.} \quad It follows from general theory (see, e.g., \cite[Chapter 1]{ShP}) that under enough general conditions an optimal stopping time in problem
(\ref{optstop}) is the first time at which the process $X_t$ exits the continuation set $\C=\{U(x)>g(x)\}$, where $\ds U(x)=\sup_{\t\in \M}\Ex g(X_\t) e^{-\r\t}$ is a value function of the problem (\ref{optstop}).

Under the above assumptions, the
domain $\C$ is an open set in $\R^1$ and, therefore, it can be
represented as a countable union of disjoint (open)
intervals. For this reason, the problem of determining
the optimal stopping time can be reduced to determining an optimal first exit time from intervals $]a,b[,\ l\le a< b\le r$,
which contain the initial point $x$ of the process $X_t$.
Necessary conditions for the optimality of such continuation intervals were obtained in \cite{Al2001}.

Frequently, it is optimal to stop when the process $X_t$ exceeds some level (threshold strategy). In this case, an optimal stopping problem can reduced to a more simple one-parametric problem of finding an optimal threshold. Similar threshold problems arise, e.g., in mathematical finance \cite{Sh}, investment models under uncertainty (real option theory) \cite{DP}, etc. Threshold strategies for optimal stopping problems in discrete time were studied in~\cite{JKS}.

However, there also exist optimal stopping problems for which the optimal strategy is not threshold (a simple example with a geometric Brownian motion
and a nondecreasing convex payoff function was given
in \cite{GuoSh}; see also \cite[Subsection 6.4]{DK}).

Another interesting question related to threshold strategies is a structure of a continuation set. In general, even if the threshold strategy is optimal in stopping problem, the continuation set $\C$ may not have a threshold structure, i.e. be not a `semi-interval' of the type $]l, p[$ (or $]p,r[$). At the same time a threshold structure of the continuation set gives, on the one hand, simple rules to stop a process starting from arbitrary initial point, and, on the other hand, allows to consider the free-boundary problem in a simple domain (semi-interval).  A conditions under which the continuation set has a threshold structure were studied in \cite{AS12}.
It derives necessary and sufficient conditions in order to the continuation set in optimal stopping problem over threshold stopping times be semi-interval. Also, this paper formulates sufficient conditions for the continuation set remains semi-interval in the problem over all stopping times. In the present paper we prove that some of the proposed in \cite{AS12} sufficient conditions are, in fact, necessary conditions.
Another sufficient conditions for the continuation set to be of threshold type were suggested in \cite{V}.

In Section 3 we give necessary and sufficient conditions for optimality in stopping problem over a class of threshold stopping times. These conditions concerning a threshold structure of optimal stopping time and continuation set are extended for the stopping problem over a class of all stopping times in Section 4. In Section 5 we demonstrate how one can use the known second-order optimality conditions for studying a relation between solutions to optimal stopping problem and to free-boundary problem.
\\

{\bf 3. Threshold stopping times.} \  In \cite{AS08}, a variational approach to solving the
optimal stopping problem for diffusion processes was
developed. In the framework of this approach, one can defined a class
of stopping times which are the first exit time out of the set (from a given family of sets), and find optimal stopping time over this class.

Let us $\tp=\inf\{t{\ge} 0: X_t\ge p\}$ be the first time when the process $X_t$ leaves the interval $]l, p[$.  We will call $\tp$ as \emph{threshold stopping time} (first exit time over threshold $p$).
Consider the optimal stopping problem (\ref{optstop}) over the
class of \emph{threshold stopping times} $\M_{\rm th}=\{\tp, \ p\in \I\}$:
\beq \label{optstop-1}
\Ex g(X_\tp) e^{-\r\tp}\to \sup_{p\in \I} .
\eeq

Let $\ds V_p(x)=\Ex g(X_\tp) e^{-\r\tp}$, and $\ds V(x)=\sup_{p\in \I} V_p(x)$ denote the optimal value in problem (\ref{optstop-1})) over threshold stopping times.

Define a \emph{continuation set over the class} $\M_{\rm th}$ as $\C_{\rm th}=\{V(x)>g(x)\}$. This definition is similar to the above definition of the continuation set for the classical optimal stopping problem (over all stopping times).

In general, even if the threshold stopping time is optimal in stopping problem (\ref{optstop}), the continuation set $\C$ may not have a threshold structure (we return to this question in Section 4).

The following theorem gives a necessary and sufficient
conditions for the optimality of  threshold stopping time in problem (\ref{optstop-1}), and also for a threshold structure of continuation set $\C_{\rm th}$.

Let the following `left-end' condition hold:
\beq\label{left-end}
\lim_{x\downarrow  l} \frac{g(x)}{\varphi(x)}=0.
\eeq

Define the function $h(p)=g(p)/\psi(p)$, where $p\in \I$.
\medskip

{\sc Theorem 1.} {\it

\emph{i)} Threshold stopping time $\t_{p^*}$ is optimal in the problem (\ref{optstop-1}) for all $x\in \I$ if and only if the following conditions hold:
\beq \label{criteria0}
h(p)\le h(p^*)\  \mbox{\rm whenever }  p<p^*;\quad h(p)\ \mbox{\rm does not increase for }  p>p^*.
\eeq

\emph{ii)} The continuation set  in problem (\ref{optstop-1}) $\C_{\rm th}$ is $]l, p^*[$, where $p^*\in\I$, if and only if the following conditions hold:
\beq \label{criteria}
h(p)<h(p^*)\  \mbox{\rm whenever }  p<p^*;\quad h(p)\ \mbox{\rm does not increase for }  p>p^*.
\eeq
}

{\it Proof.} \ Let us consider $\t_{a,p}=\inf\{t{\ge} 0: X_t\notin ]a,p[\}$, where $l<a<p$.
The known formulas for a solution to two-sided boundary problem (see, e.g. \cite[Lemma 4.3]{DK}, or \cite[Section 3]{Al2001}) give the following representation:
\begin{eqnarray*}
&&V_{a,p}(x):= \Ex g(X_{\t_{a,p}}) e^{-\r\t_{a,p}}=g(a)u_1(x,a,p)+ g(p) u_2(x,a,p),\quad \mbox{\rm where}\\[6pt]
&&u_1(x,a,p)=\frac{\psi(x)\varphi(p)- \psi(p)\varphi(x)} {\psi(a)\varphi(p)- \psi(p)\varphi(a)},\quad u_2(x,a,p)=\frac{\psi(a)\varphi(x)- \psi(x)\varphi(a)} {\psi(a)\varphi(p)- \psi(p)\varphi(a)}.
\end{eqnarray*}

Thus, letting $a\downarrow l$,  the `left-end' condition (\ref{left-end}) implies the following formula:
\beq\label{representation}
V_p(x)=\left\{
\begin{array}{ll}
h(p)\psi(x), & \mbox{\rm for } x<p,\\[4pt]
g(x), &    \mbox{\rm for } x\ge p.
\end{array}
\right.
\eeq

i) Let  (\ref{criteria0}) hold. Take arbitrary $x, p\in \I$.

If $x<\min(p,p^*)$ then $V_p(x)=h(p)\psi(x)\le h(p^*)\psi(x)=V_{p^*}(x)$.

If $x\ge \max(p,p^*)$ then $V_p(x)=g(x)=V_{p^*}(x)$.

If $p\le x<p^*$ then $V_p(x)=g(x)=h(x)\psi(x)\le h(p^*)\psi(x)=V_{p^*}(x)$.

If $p^*\le x<p$ then $V_p(x)=h(p)\psi(x)\le h(x)\psi(x)=g(x)=V_{p^*}(x)$.

Now, let $V_p(x)\le V_{p^*}(x)$ for all $x$ and $p$. Then for $p<p^*$ and $x<p$ we have: $V_p(x)=h(p)\psi(x)\le V_{p^*}(x)=h(p^*)\psi(x)$, i.e. $h(p)\le h(p^*)$. If $p^*<p_1<p_2$ then $V_{p_2}(p_1)=h(p_2)\psi(p_1)\le V_{p^*}(p_1)=g(p_1)=h(p_1)\psi(p_1)$, therefore, $h(p_2)\le h(p_1)$.

ii) The representation (\ref{representation}) implies
\begin{eqnarray*}
V(x)&=&\sup_p V_p(x)=\max\{\sup\limits_{p\le x} V_p(x), \sup\limits_{ p> x} V_p(x)\} \\
&=&\max\{g(x), \psi(x)\cdot\sup\limits_{ p> x} h(p)\}=\psi(x)\max\{h(x), \sup\limits_{ p> x} h(p)\}.
\end{eqnarray*}

Therefore,  the continuation set $\C_{\rm th}$ and the stopping set $\S_{\rm th}$ in problem (\ref{optstop-1}) can be written as follows:
$$
\S_{\rm th}=\{x:\ V(x)=g(x)\}=\{x:\ \sup\limits_{ p> x} h(p)=h(x)\},
$$
$$
\C_{\rm th}=\{x:\ V(x)>g(x)\}=\{x:\ \sup\limits_{ p> x} h(p)>h(x)\}.
$$

Conditions (\ref{criteria}) imply that $h(p^*)\ge h(p)$ for $p>p^*$. Therefore, for any $x<p^*$ we have $\sup\limits_{ p> x} h(p)\ge h(p^*)>h(x)$, i.e. $x\in \C_{\rm th}$. If $x\ge p^*$, then $\sup\limits_{ p> x} h(p)= h(x)$ and, hence, $x\in \S_{\rm th}$. It means that $\C_{\rm th}=]l,p^*[$.

Conversely, let $\C_{\rm th}=]l, p^*[$, and, therefore, $\S_{\rm th}=[p^*,r[$.

Suppose that $h(x_1)<h(x_2)$ for some $p^*<x_1<x_2$. Then $\sup\limits_{ p> x_1} h(p)\ge h(x_2)>h(x_1)$, that contradicts to $x_1\in \S_{\rm th}$. Therefore, $h(p)$ decreases for $p>p^*$.

Suppose that $h(x)\ge h(p^*)$ for some $x<p^*$. Denote $\ds\bar x=\min\{y\in [x,p^*]: \ h(y)=\max_{x\le p\le p^*} h(p)\}$.
If $\bar x<p^*$, then $\sup\limits_{ p> \bar x} h(p)=\sup\limits_{ [\bar x, p^*]} h(p) = h(\bar x)$, that contradicts to $\bar x\in \C_{\rm th}$. If $\bar x=p^*$, then $h(x)<h(p^*)$, but it contradicts to the assumption $h(x)\ge h(p^*)$. Hence, $h(x)< h(p^*)$ for all $x<p^*$.

Theorem proved.

\medskip
So, the optimal threshold $p^*$ is a point of maximum for the function $h(p)$. This implies the necessity (under minor assumptions) of
the well-known smooth pasting principle.
\medskip

{\bf Corollary.} {\it Suppose that $\t_{p^*}$, where $p^*\in \I$, is the optimal stopping time in the
problem (\ref{optstop-1}) and the function $g(x)$ is differentiable at the point $p^*$. Then the function $V(x)$ is differentiable at the
point $p^*$, and ${V}'(p^*){=}g'(p^*)$.
}\medskip

{\bf Remark.} If the function  $g(x)$ has only one-sided derivatives $g'(p^*{-}0)$ and $g'(p^*{+}0)$, then instead of ``smooth pasting"\ the following inequalities hold:
$$ g'(p^*{+}0)=V'(p^*{+}0)\le V'(p^*{-}0))\le g'(p^*{-}0).$$

{\it Proof. } Since $p^*$ is a point of maximum for the function $h$, then
$$
h'(p^*{-}0)\ge 0\ge h'(p^*{+}0).
$$
Therefore, we have
$$
V'(p^*{-}0)=h(p^*)\psi'(p^*)=g'(p^*{-}0)-h'(p^*{-}0)\psi(p^*)\le g'(p^*-0),
$$
$$
V'(p^*{+}0)= g'(p^*{+}0)=h'(p^*{+}0)\psi(p^*)+h(p^*)\psi'(p^*)\le h(p^*)\psi'(p^*)=V'(p^*{-}0).
$$
\medskip

The necessity of the smooth pasting condition under
some additional constraints on the process was shown
in \cite{ShP}. A result similar to ours was obtained in \cite{V}.
\medskip

{\bf 4. Threshold structure of optimal stopping time and continuation set.}\
In this Section we give conditions under which an optimal stopping time and continuation set  have a  threshold structure in general problem (\ref{optstop}) over all stopping times.
Recall that a continuation set in problem (\ref{optstop}) is  $\C=\{x:\,U(x)>g(x)\}$, where $\ds U(x)=\sup_{\t\in \M}\Ex g(X_\t) e^{-\r\t}$, and $\S=\{x:\,U(x)=g(x)\}$ is a stopping set.

\medskip

{\sc Theorem 2.}

{\it \emph{i)} Let for some $p^*\in \I$ a payoff function $g$ be twice continuously differentiable on interval $]p^*, r[$ and there exists $g'(p^*{+}0)$. If the following conditions hold:
\begin{eqnarray}
&&\hskip-1cm h(p)\le h(p^*) \quad\mbox{\rm  for }  p<p^*;\label{criteria01}\\[3pt]
&&\hskip-1cm {\psi'(p^*)}g(p^*)\ge {\psi(p^*)}g'(p^*{+}0);\label{criteria011}\\[3pt]
&& \hskip-1cm
 \L g(p)\le \r g(p) \quad \mbox{\rm  for } p>p^*, \label{criteria2}
\end{eqnarray}
then $\t_{p^*}$ is the optimal stopping time in problem (\ref{optstop}) for all $x\in \I$.\smallskip

 \emph{ii)} If \ $\t_{p^*}$ is the optimal stopping time in problem (\ref{optstop}) for all $x\in \I$, $g$ be twice continuously differentiable on interval $]p^*, r[$, and there exists $g'(p^*{+}0)$, then conditions (\ref{criteria01})--(\ref{criteria2}) hold.
}
\medskip

If inequality in (\ref{criteria01}) is strict, one can get a similar results  for a continuation set  $\C$.\medskip

{\sc Theorem 3.}

{\it \emph{i)} Let for some $p^*\in \I$  a payoff function $g$ be twice continuously differentiable on interval $]p^*, r[$, and there exists $g'(p^*{+}0)$. If
\beq
h(p)<h(p^*)\  \mbox{\rm whenever } p<p^*,\label{criteria1}
\eeq
and (\ref{criteria011})--(\ref{criteria2}) hold, then $\C=]l,p^*[$.\smallskip

 \emph{ii)} Let $\limsup_{x\uparrow r}\max(0,h(x))=0$. If \ $\C=]l,p^*[$, a payoff function $g$ be twice continuously differentiable on interval $]p^*, r[$, and there exists $g'(p^*{+}0)$, then conditions (\ref{criteria1}) and (\ref{criteria011})--(\ref{criteria2}) hold.
}
\medskip

{\it Proof of Theorem 2.} \ i) Let (\ref{criteria01})--(\ref{criteria2}) hold. Take the function
$$
\Phi(x)=V_{p^*}(x)=\left\{
\begin{array}{ll}
h(p^*)\psi(x), & \mbox{\rm for } x<p^*,\\[4pt]
g(x), &    \mbox{\rm for } x\ge p^*.
\end{array}
\right.
$$

Obviously, $ U(x)\ge \Phi(x)$.

On the other hand, (\ref{criteria01}) implies $h(p^*)\psi(x)\ge h(x)\psi(x)=g(x)$, therefore $\Phi(x)\ge g(x)$ for all $x$. This inequality and Ito--Tanaka formula (see, e.g. \cite{KS}) imply: for any stopping time $\t\in \M$ and $N>0$, $\t_N=\t\wedge N$
\begin{eqnarray}
\Ex g(X_{\t_N}) e^{-\r{\t_N}}&\le& \Ex \Phi(X_{\t_N}) e^{-\r{\t_N}}=\Phi (x)+\Ex \int_0^{\t_N}(\L\Phi -\r\Phi)(X_t)e^{-\r t}dt\nonumber\\
&+&\frac12 \s^2(p^*)[\Phi'(p^*{+}0)-\Phi'(p^*{-}0)]\Ex \int_0^{\t_N} e^{-\r t}dL_t(p^*),
\end{eqnarray}
where $L_t(p^*)$ is the local time at $p^*$.

Define $T_1=\{0\le t\le {\t_N}:\, X_t < p^*\}$, $T_2=\{0\le t\le {\t_N}:\, X_t > p^*\}$. We have: $\L\Phi(X_t) -\r\Phi(X_t) =h(p^*) [\L\psi(X_t)-\r\psi(X_t)]=0$ for $t\in T_1$, and $\L\Phi(X_t) -\r\Phi(X_t) =\L g(X_t) -\r g(X_t)\le 0$ for $t\in T_2$ (see (\ref{criteria2})).

By definition we have:
$$
\Phi'(p^*{+}0)=g'(p^*{+}0), \quad \Phi'(p^*{-}0)= h(p^*)\phi'(p^*).
$$
Thus, due to (\ref{criteria011}) $\Phi'(p^*{+}0)-\Phi'(p^*{-}0)\le 0$.

Then
\begin{eqnarray*}
&&\Ex g(X_{\t_N}) e^{-\r\t_N}\le \Phi (x)+\Ex \left( \int\limits_{T_1}(\L\Phi -\r\Phi)(X_t)e^{-\r t}dt \right. \\
&&+ \left.\int\limits_{T_2}(\L\Phi -\r\Phi)(X_t)e^{-\r t}dt\right)+(...)[\Phi'(p^*{+}0)-\Phi'(p^*{-}0)] \le \Phi (x).
\end{eqnarray*}
Letting $N\to \infty$ we get $\ds \Ex g(X_{\t}) e^{-\r\t}\le \Phi (x)$.
It follows from this inequality, that $ U(x)\le \Phi(x)$.

Therefore, $ U(x)= \Phi(x)=V_{p^*}(x)$, i.e. $\t_{p^*}$ is the optimal stopping time in problem (\ref{optstop}) for all $x$.

ii) Let us note, that $\t_{p^*}$ will be an optimal stopping time in the problem (\ref{optstop-1}) also. Thus, (\ref{criteria01}) follows immediately from Theorem 1.

Assume, that inequality (\ref{criteria2}) is not true at some point $x_0 > p^*$, i.e. $\L g(x)> \r g(x)$ in some interval $J\subset ]p^*,r[$ (by virtue of continuity). For some $\tilde x\in J$ and $N>0$ let us define $ \t=\inf\{t\ge 0:\, X_t\notin J\}$, where process $X_t$ starts from the point  $\tilde x$, and $\t_N=\t\wedge N$. Then
$$
\E^{\tilde x} g(X_{\t_N}) e^{-\r \t_N}= g (\tilde x)+\E^{\tilde x} \int_0^{\t_N}(\L g -\r g)(X_t)e^{-\r t}dt >g (\tilde x).
$$
Therefore, $U(\tilde x) {>} g(\tilde x)$ that contradicts to $U(\tilde x){=}V_{p^*}(\tilde x){=}g(\tilde x)$ (since $\tilde x>p^*$).

This completes the proof.\\

{\it Proof of Theorem 3.} \ i) As in the proof of Theorem 2 we get  $ U(x)= \Phi(x)=V_{p^*}(x)$, and, therefore,
\begin{eqnarray*}
\C&=&\{x:\,U(x)>g(x)\}=\{x:\,x<p^*,\ h(p^*)\psi(x)>g(x)\}\\
&=& \{x:\,x<p^*,\ h(p^*)>h(x)\}=\{x:\,x<p^*\}
\end{eqnarray*}
due to (\ref{criteria1}).

ii) Proposition 5.7 in \cite{DK} implies that $\t_{p^*}$ is the optimal stopping time in the problem (\ref{optstop}), and, therefore, in the problem (\ref{optstop-1}). Thus, the value functions in these problems ($U(x)$ and $V(x)$) are the same, and, therefore, $\C = \C_{\rm th}$. Now, (\ref{criteria1}) follows immediately from Theorem 1.

If (\ref{criteria2}) is not true at some point $x_0 > p^*$, then repeating arguments as above in the proof of Theorem 2 one can get that $U(\tilde x) {>} g(\tilde x)$ for some $\tilde x>p^*$. But this contradicts to $\tilde x \in\S$.

The proof is complete.

\medskip
{\bf Corollary.} {\it
For the case of linear payoff function $g(x) = x - c$ threshold stopping time $\t_{p^*}$ is optimal in the stopping problem (\ref{optstop}) if and only if the following conditions hold:
\begin{eqnarray}
&&
\frac{p-c}{\psi(p)}\le \frac{p^*-c}{\psi(p^*)}\quad \mbox{\rm  for } p<p^*;\label{criteria110}\\[3pt]
&& \psi'(p^*)(p^*-c) \ge \psi(p^*);
\label{criteria11}\\
&& a(p)\le \r (p- c) \quad \mbox{\rm  for } p>p^*. \label{criteria21}
\end{eqnarray}
where $a(p)$ is the drift function of process $X_t$.

Conditions (\ref{criteria110}) with strict inequality and (\ref{criteria11})--(\ref{criteria21}) will be necessary and sufficient in order to $\C=]l,p^*[$.
}
\medskip

Note that, in the classical case of real option ($g(x) = x - c$) or American call option ($g(x) = \max(0,x - c)$), when
the process $X_t$ is a geometric Brownian motion, the conditions (\ref{criteria110})--(\ref{criteria21}) hold automatically. Although the function  $(x - c)_+$ is not smooth at point $x=c$, condition (\ref{criteria11}) implies $p^*>c$, so that it is twice differentiable for $p>p^*$.\\

\textbf{Remark.} Maximality of $h(p)$ at $p^*$ and condition (\ref{criteria2}) imply that $h(p)\downarrow$ for $p>p^*$.
\smallskip

{\it Proof.} Indeed, for $p^* < x<y$ and process  $X_t$ starting from the point  $x$, define $\t_N=\inf\{t\ge 0:\, X_t\notin ]p^*,y[\}\wedge N$ and function $f(z)=g(z)-h(y)\psi(z)$ for $z\in \I$. Obviously, $\L f(z)-\r f(z)=\L g(z)-\r g(z)$. By Dynkin's formula and (\ref{criteria2}) we have
\begin{eqnarray}
\E^{x} f(X_{\t_N}) e^{-\r \t_N}&=& f (x)+\E^{x} \int_0^{\t_N}(\L f -\r f)(X_t)e^{-\r t}dt \nonumber\\
&=& f (x)+\E^{x} \int_0^{\t_N}(\L g -\r g)(X_t)e^{-\r t}dt\nonumber\\
&\le& f(x)=\psi(x)[h(x)-h(y)].\label{12}
\end{eqnarray}

On the other hand, if $\t_1=\inf\{t\ge 0:\, X_t=p^*\}$, $\t_2=\inf\{t\ge 0:\, X_t=y\}$, then
\begin{eqnarray*}
\E^{x} f(X_{\t_N}) e^{-\r \t_N}&=& \E^{x} f(X_{N}) e^{-\r N} (x)\chi\{\t_1,\t_2>N\}\\
&+& \E^{x} f(p^*) e^{-\r \t_1} (x)\chi\{\t_1 <N,\t_2>\t_1\}\\
&+& \E^{x} f(y) e^{-\r \t_2} (x)\chi\{\t_2 <N,\t_1>\t_2\}=E_1+E_2+E_3,
\end{eqnarray*}
where $\chi\{\cdot\}$ is the indicator of event $\{\cdot\}$. Further,\\
$\ds E_1=e^{-\r N}\E^{x} f(X_{N})\chi\{\t_1,\t_2>N\}\ge -e^{-\r N} \max_{p^*\le z\le y}f(z)$,\\
$\ds E_2{=}f(p^*)\E^{x}e^{-\r \t_1} (x)\chi\{\t_1 {<}N,\t_2{>}\t_1\}{\ge} 0$ (since $f(p^*){=}\psi(p^*)[h(p^*){-}h(y)]{\ge} 0$),\\
$\ds E_3{=}f(y)\E^{x}e^{-\r \t_2} (x)\chi\{\t_2 {<}N,\t_1{>}\t_2\}=0$.

Thus, (\ref{12}) implies
$$
\psi(x)[h(x)-h(y)]\ge -e^{-\r N} \max_{p^*\le z\le y}f(z),
$$
and taking $N\to\infty$ we get $h(x)\ge h(y)$. \quad $\Box$  \\

{\bf 5. Free-boundary problem and second-order conditions.} \  The solution to the optimal stopping problem (\ref{optstop}),
i.e. the continuation set and the value function, is usually sought as a solution to the free-boundary problem for the differential operator $\L$ (FB problem). The obtained solution is considered as a ``candidate'' for the solution to the optimal stopping problem, which needs an additional verification for optimality.

As applied to threshold strategies considered in this
paper, free-boundary problem looks as follows: to find a threshold $p^*$ (determining the continuation set) and a function $U(x),\ l<x<p^*$, such that
\begin{eqnarray}
&&\L U(x)=\r U(x), \quad l<x< p^*;\label{Stefan1} \\
&& U( p^*-0)= g(p^*), \label{Stefan2}\\
&& U'(p^*-0)=g'(p^*).\label{Stefan3}
\end{eqnarray}

Conditions (\ref{Stefan1}) and (\ref{Stefan2}) are satisfied for the function
$U(x)=V_{p^*}(x)= h(p^*)\psi(x)$, \  $x <p^*$. The smooth-pasting condition (\ref{Stefan3}) at the point $p^*$ is equivalent to the stationarity of the function $h(p)$ at this point, i.e. to $h'(p^*)=0$.

On the other hand, it follows from the above results
that the threshold $p^*$, which determines the boundary of continuation set, is a point of maximum for the function $h(p)$.

These remarks make it possible to use classical extremum conditions for analyzing the relationship between solutions to free-boundary problem and to optimal stopping problem.

Example 1 below demonstrates how a difference between stationarity and maximality may be applied to the problems under consideration.
\medskip

{\bf Example 1} (a solution to free-boundary problem \emph{may not
be} a solution to the optimal stopping problem, see \cite{AS08}).

Consider the geometric Brownian motion $dX_t=X_t(0.5 dt+ dw_t)$ with $X_0=x$, the payoff function $g(x){=}
(x-1)^3+x^\d$ \ ($\d>0,\ x\ge 0$), and the discount
 $\r=\d^2/2$. Note that the function $g$ is smooth and monotonically increases for all $\d>0$, and $\psi(x)=x^\d$.

For $\d< 3$, FB problem (\ref{Stefan1})--(\ref{Stefan3}) has the unique solution $ U(x)=V_1(x)=x^\d$,\  $p^*=1$. However $\t_1$ is not an optimal stopping time, because $\ds V_p(x)=[(p-1)^3p^{-\d}+1]x^\d\to \infty$ as $p\to \infty$ for any $x>0$. Thus,a solution to
free-boundary problem does exist, while the optimal stopping problem has no solution.

For $\d= 3$, FB problem has the unique solution
$ U(x)=V_1(x)=x^3$,\  $p^*=1$, which cannot be a solution
to the optimal stopping problem, because $V_p(x)\uparrow V(x)=2x^3$  as $p\to \infty$. Thus, the solution $U(x)$ to free-boundary problem is strictly less than the value function $V(x)$ in problem (\ref{optstop-1}), which, although finite, is not attained by any threshold strategy.

For $\d> 3$, FB problem (\ref{Stefan1})--(\ref{Stefan3}) has two solutions, $
U(x){=}V_1(x){=} x^\d$, \ $p^*=1$, and $
U(x){=}V_{p_\d}(x){=}h(p_\d)x^\d$, \ $p^*{=}p_\d{=}\d/(\d{-}3)$. Note that $V_{p_\d}(x) > V_1(x)$ for $0<x<p_\d$. Thus, one of the solutions to FB problem (with boundary $p^*=1$) does not give a solution
to the optimal stopping problem (which exists, unlike
in the preceding situations). It can be shown (by applying, e.g., Theorem 2) that the first exit time over the boundary $p_\d$ is optimal among all stopping times $\M$.

\medskip

\medskip

As it was shown above (see (\ref{representation}))  a solution to free-boundary problem (\ref{Stefan1})--(\ref{Stefan3}) in one-dimensional case (under the condition (\ref{left-end})) has the following simple structure:
\beq
U(x)=V_{p^*}(x)= h(p^*)\psi(x),\quad x<p^*, \label{Stefan-4}
\eeq
where $h'(p^*)=0$.

Formula (\ref{Stefan-4}) implies: $U''(p^*-0)=h(p^*)\psi''(p^*)$.
Since $g(p)=h(p)\psi(p)$ and $h'(p^*)=0$, then
$g''(p^*)=h''(p^*)\psi(p^*)+h(p^*)\psi''(p^*)$. Therefore,
\beq
\psi(p^*)h''(p^*)=g''(p^*)-U''(p^*-0).\label{Stefan-2}
\eeq
It means that $\mbox{\rm sign\,} h''(p^*)=\mbox{\rm sign\,}[g''(p^*)-U''(p^*-0)]$.

On the other hand, if $x<p^*$ then condition $h'(p^*)=0$ and a sign of $h''(p^*)$ are characterized a local extremum (in $p$) of the function $V_p(x)$ at the point $p^*$. Hence, applying known second-order optimality conditions one can derive  the following

\medskip

{\bf Proposition.} {\it Suppose that a pair $(U(x),p^*)$ is a solution to free-boundary problem (\ref{Stefan1})-(\ref{Stefan3}) and $g\in C^2(O(p^*))$ in some neighborhood of $p^*$. Then:

\emph{(i)} if \ $V(x) = U(x)\chi\{x < p^*\}+g(x)\chi\{x \ge p^*\}$  is the value function in the optimal stopping problem (\ref{optstop-1}), then $U''(p^*-0) \ge g''(p^*)$;

\emph{(ii)} if \ $U''(p^*-0)>g''(p^*)$ and $x < p^*$, then $p^*$ is a point of strict local maximum (in $p$) of the function $V_p(x)$;

\emph{(iii)} if \ $U''(p^*-0)<g''(p^*)$ and $x < p^*$, then $p^*$ is a point of strict local minimum (in $p$) of the function $V_p(x)$.
}\\

\textbf{Remark.} The case $U''(p^*-0)=g''(p^*)$ needs the additional considerations and maybe using a high-order conditions.\\

If free-boundary problem has several solutions, then the
above proposition makes it possible to discard those
solutions which surely cannot solve the optimal stopping problem (e.g., with boundaries which are the points of local minimum
of the function  $V_p(x)$ in $p$).
\medskip

{\bf Example 2.} Consider the process of geometric
Brownian motion $dX_t=X_t(\a dt+\s dw_t),\ X_0=x$, the discount $\r=\s^2+2\a$, and the smooth monotonically increasing payoff function defined by
$$
g(x)=\left\{
\begin{array}{ll}
[(x-1)^2+1]x^2, & \mbox{ for }\ 0\le x\le 2,\\[4pt]
x-9+\frac{15}{4} x^2, & \mbox{ for }\  x > 2.
\end{array}
\right.
$$

In this case, free-boundary problem has two solutions:

(a) $ U_1(x)=x^2,\ p_1^*=1$, and

(b) $U_2(x)=\frac{34}9 x^2,\ p_2^*=18$.
\smallskip

The second-order condition gives: \  $U_1''(p_1^*-0)=2 < g''(p_1^*)=4$ \  and \  $U_2''(p_2^*-0)=\frac{68}9 > g''(p_2^*)=7.5$. Therefore, solution (a) surely cannot be a solution to the optimal stopping problem (even over the class of threshold strategies), while (b) is indeed the solution
to the optimal stopping problem, which can easily
shown by using Theorem 2.
\\

{\bf Acknowledgement.} The authors are thankful to Ernst Presman for helpful remarks and discussions.

\end{document}